\documentclass[12pt,twoside]{article} 
\usepackage{amsmath}
\usepackage{amssymb}

\newcommand{\Q}{\mbox{{\boldmath $Q$}}}

\title{%
\Large\bfseries
Some identities of symmetry for the generalized Bernoulli numbers and polynomials
}
\author{%
\vspace{5mm}
By \\
Taekyun {\sc Kim} } \vspace{15mm}
\def\Q{\mbox{\boldmath $Q$}}

\date{%
\begin{minipage}{14cm}%
\vspace{10mm}
\normalsize
{\bfseries {\small Abstract.}} {\small
In this paper, by the properties of $p$-adic invariant integral on $\mathbb{Z}_{p}$,
we establish various identities concerning the generalized Bernoulli numbers and polynomials.
From the symmetric properties of $p$-adic invariant integral on $\mathbb{Z}_{p}$,
 we give some interesting relationship between the power sums and the generalized Bernoulli polynomials.
} 
\\
\end{minipage} \\
\begin{minipage}{18cm}
\normalsize
{\bfseries {\small 2000 Mathematics Subject Classification:}}
{\small
11B68, 11M38, 11S80.
 }
\noindent \\
{\bfseries {\small Key Words and Phrases:}}
{\small $p$-adic invariant integral, Bernoulli numbers, Bernoulli\\
polynomials.
}
\\
\end{minipage}
}
\newcommand{\SECTION}[2]{%
  \vspace{5mm}\par\noindent $\S$#1. {\bf #2}~}

\newcommand{\THM}[1]{%
  \par\vspace{5mm}\par\noindent {\scshape #1.~}}

\newcommand{\LEM}[1]{%
    \par\vspace{5mm}\par\noindent {\scshape #1.}}

\newcommand{\REM}[1]{%
    \par\vspace{5mm}\par\noindent {\scshape #1.}}

\newcommand{\pn}{\par\noindent}

\newcommand{\mpn}{\medskip\par\noindent}
\newcommand{\bpn}{\bigskip\par\noindent}
\begin{document}
\maketitle \pagestyle{myheadings} \markboth{{\footnotesize {\hfill
\rm Some identities of symmetry for the generalized Bernoulli
numbers and polynomials \hfill}}} {{ \small {\hfill \rm Taekyun {\sc
Kim}
  \hfill} }}
 \vspace* {-10mm}
 \vspace* {-25mm}
 \vspace* {-55mm}

\newcommand\QF[1]{\Q \hspace{-1mm} \left( \hspace{-1mm}\sqrt{#1} \right) }
\newenvironment{THEOREM}[1]{%
    \vspace{5mm}\par\noindent{\bf Theorem. #1}~}{}
\newenvironment{REMARK}[1]{%
    \vspace{5mm}\par\noindent{\bf Remark #1}~}{}
\newenvironment{LEMMA}[1]{%
    \vspace{5mm}\par\noindent{\bf Lemma #1}~}{}
\newenvironment{PROPOSITION}[1]{%
    \vspace{5mm}\par\noindent{\bf Proposition #1}~}{}
\newenvironment{COROLLARY}[1]{%
    \vspace{5mm}\par\noindent{\bf Corollary. #1}~}{}
\newenvironment{DEFINITION}[1]{%
    \vspace{5mm}\par\noindent{\bf Definition #1}~}{}
\newenvironment{PROOF.}{%
    \vspace{5mm}\par\noindent{\bf Proof.}~}{%
    \hfill\hbox{\rule[-2pt]{3pt}{6pt}}%
    \par\vspace{5mm}}
\newenvironment{EXAMPLES}{%
    \vspace{5mm}\par\noindent{\sc Examples}.~}{}
\newenvironment{ACKNOWLEDGEMENTS}{%
    \vspace{5mm}\par\noindent{\bf Acknowledgements.}~}{}
\newcommand{\Mod}[1]{\,(\text{\mbox{\rm mod}}\;#1)}
\newcommand{\NUM}[1]{{[#1]}}
\baselineskip 6mm
\thispagestyle{empty}
\vspace{2cm}
\vspace{5cm}
\medskip
\par
\begin{center}\SECTION{1}{ Introduction}
\end{center}
\vspace{3mm}
\par
Let $p$ be a fixed prime number. Throughout this paper, the symbols
$\mathbb{Z}, \mathbb{Z}_p, \mathbb{Q}_p$, and $\mathbb{C}_p$ will
denote the ring of rational integers, the ring of $p$-adic integers,
the field of $p$-adic rational numbers, and the completion of
algebraic closure of $\mathbb{Q}_p $, respectively. Let $\mathbb{N}$
be the set of natural numbers and $\mathbb{Z}_{+}=\mathbb{N}\cup \{
0 \}$. Let $v_{p}$ be the normalized exponential valuation of
$\mathbb{C}_{p}$ with $|p|_{p}=p^{-v_{p}(p)}= 1/p$. Let
$UD(\mathbb{Z}_{p})$ be the space of uniformly differentiable
function on $\mathbb{Z}_{p}$. For $f\in UD(\mathbb{Z}_{p})$, the
$p$-adic invariant integral on $\mathbb{Z}_{p}$ is defined as
\begin{equation*}
\tag{1}
I(f)=\int_{\mathbb{Z}_{p}} f(x)dx= \lim_{N \rightarrow \infty}
\frac{1}{p^{N}} \sum_{x=0}^{p^{N} -1} f(x), \quad (\hbox{see [6]}).
\end{equation*}
From the definition (1), we have
\begin{equation*}
\tag{2}
I_{1}(f_{1})=I_{1}(f)+ f^{\prime}(0), \hbox{ where \ $f^{\prime}(0)=\frac{df(x)}{dx}|_{x=0}$ \ and \                 $f_{1}(x)=f(x+1)$}.
\end{equation*}
Let $f_{n}(x)=f(x+n), \ (n \in \mathbb{N})$. Then we can derive the following equation (3) from (2).
\begin{equation*}
\tag{3}
I(f_{n})=I(f)+\sum_{i=0}^{n} f^{\prime}(i), \quad (\hbox{see [6]}).
\end{equation*}
It is well known that the ordinary Bernoulli polynomials $B_{n}(x)$ are defined as
\begin{equation*}
\frac{t}{e^{t}-1}e^{xt}=
\sum_{n=0}^{\infty}B_{n}(x)\frac{t^{n}}{n!}, \quad \hbox {( see [1-25] )} ,
\end{equation*}
and the Bernoulli number $B_{n}$ are defined as $B_{n}=B_{n}(0)$.
\par
Let $d$ a fixed positive integer. For $n \in \mathbb{N}$, we set
$$
\displaystyle{X=X_{d}=\lim_{\overleftarrow{N}}\Big( \mathbb{Z}/dp^{N}\mathbb{Z}\Big)},
\quad
X_{1}=\mathbb{Z}_{p};
$$
$$
X^{*}=\bigcup_{0<a<dp , \atop (a,p)=1} (a+dp\mathbb{Z}_{p});
$$
$$
a+dp^{N}\mathbb{Z}_{p}=\{ \ x \in X |\ x\equiv a\Mod{dp^{N}} \ \},
$$
where $a \in \mathbb{Z}$ lies in $0 \leq a <dp^{N}$.
In [6], it is known that
\begin{equation*}
\int_{X}f(x)dx=\int_{\mathbb{Z}_{p}}f(x)dx, \quad  \hbox{ for $f\in UD(\mathbb{Z}_{p})$}.
\end{equation*}
Let us take $f(x)=e^{tx}$. Then we have
\begin{equation*}
\int_{\mathbb{Z}_{p}}e^{tx}dx=\frac{t}{ e^{t}-1}=\sum_{n=0}^{\infty}B_{n}\frac{t^{n}}{n!}.
\end{equation*}
Thus, we note that
\begin{equation*}
\int_{\mathbb{Z}_{p}}x^{n}dx=B_{n}, \quad \ \hbox {$n \in \mathbb{Z}_{+}$,} \quad (\hbox{see [1-25]}).
\end{equation*}
Let $\chi$ be the Dirichlet's character with conductor $d\in \mathbb{N}$. Then the generalized Bernoulli polynomials attached to $\chi$ are defined as
\begin{equation*}
\tag{4} \sum_{a=1}^{d}\frac{\chi(a)te^{at}}{e^{dt}-1}e^{xt}
=\sum_{n=0}^{\infty}B_{n,\chi}(x)\frac{t^{n}}{n!},  \quad \hbox {(
see [22] )},
\end{equation*}
and the generalized Bernoulli numbers attached to $\chi, \ B_{n,\chi}$ are defined as $B_{n,\chi}=B_{n,\chi}(0)$.
\par
In this paper, we investigate the interesting identities of symmetry for the generalized Bernoulli numbers and polynomials attached to $\chi$ by using the properties of $p$-adic invariant integral on $\mathbb{Z}_{p}$. Finally, we will give relationship between the power sum polynomials and the generalized Bernoulli numbers attached to $\chi$.
\newpage
\begin{center}
\SECTION{2}{Symmetry of power sum and the generalized Bernoulli
polynomials }
\end{center}
\vspace{5mm}
\par
Let $\chi$ be the Dirichlet character with conductor $d\in \mathbb{N}$.
From (3), we note that
\begin{equation*}
\tag{5}
\int_{X}\chi(x)e^{xt}dx=
\frac{t\sum_{i=0}^{d-1}\chi(i)e^{it}}{e^{dt}-1}=
\sum_{n=0}^{\infty}B_{n,\chi}\frac{t^{n}}{n!},
\end{equation*}
where $B_{n,\chi}(x)$ are $n$-th generalized Bernoulli numbers attached to $\chi$. Now, we also see that the generalized Bernoulli polynomials attached to $\chi$ are given by
\begin{equation*}
\tag{6} \int_{X}\chi(y)e^{(x+y)t}dy=
\frac{\sum_{i=0}^{d-1}\chi(i)e^{it}}{e^{dt}-1}e^{xt}=
\sum_{n=0}^{\infty}B_{n,\chi}(x)\frac{t^{n}}{n!}.
\end{equation*}
By (5) and (6), we easily see that
\begin{equation*}
\tag{7}
\int_{X}\chi(x)x^{n}dx=B_{n,\chi}, \quad \hbox{and} \quad
\int_{X}\chi(y)(x+y)^{n}dy=B_{n,\chi}(x).
\end{equation*}
From (6), we have
\begin{equation*}
\tag{8}
B_{n,\chi}(x)=\sum_{\ell=0}^{n}\binom{n}{\ell}B_{\ell,\chi}x^{n-\ell}.
\end{equation*}
From (6), we can also derive
\begin{equation*}
\int_{X}\chi(x)e^{xt}dx=\sum_{i=0}^{d-1}\chi(i)\frac{t}{e^{dt}-1}e^{(\frac{i}{d})dt}
=\sum_{n=0}^{\infty}\Big(d^{n}\sum_{i=0}^{d-1}\chi(i)B_{n}(\frac{i}{d})\Big)\frac{t^{n}}{n!}.
\end{equation*}
Therefore, we obtain the following lemma.
{\LEM {Lemma{1}}} {\it For $n\in \mathbb{Z}_{+}$, we have
\begin{equation*}
\int_{X}\chi(x)x^{n}dx=B_{n,\chi}=d^{n}\sum_{i=0}^{d-1}\chi(i)B_{i}\Big(\frac{i}{d}\Big).
\end{equation*}
\rm}
\medskip
\par
We observe that
\begin{equation*}
\tag{9}
\frac{1}{t}\Big(\int_{X}\chi(x)e^{(nd+x)t}dx-\int_{X}e^{xt}\chi(x)dx \Big)
=\frac{nd\int_{X}\chi(x)e^{xt}dx}{\int_{X}e^{ndxt}dx}
=\frac{e^{ndt}-1}{e^{dt}-1}\Big( \sum_{i=0}^{d-1}\chi(i)e^{it}\Big).
\end{equation*}
Thus, we have
\begin{equation*}
\tag{10}
\frac{1}{t}\Big(\int_{X}\chi(x)e^{(nd+x)t}dx-\int_{X}e^{xt}dx \Big)
=\sum_{k=0}^{\infty}\Big( \sum_{\ell=0}^{nd-1}\chi(\ell)\ell^{k}\Big)\frac{t^{k}}{k!}.
\end{equation*}
Let us define the $p$-adic functional $T_{k}(\chi, n)$ as follows:
\begin{equation*}
\tag{11}
T_{k}(\chi,n)=\sum_{\ell=0}^{n}\chi(\ell)\ell^{k}, \quad \hbox{for \ $k\in \mathbb{Z}_{+}$}.
\end{equation*}
By (10) and (11), we see that
\begin{equation*}
\tag{12}
\frac{1}{t}\Big(\int_{X}\chi(x)e^{(nd+x)t}dx-\int_{X}e^{xt}dx \Big)
=\sum_{n=0}^{\infty}\Big( T_{k}(\chi,nd-1)\Big)\frac{t^{k}}{k!}.
\end{equation*}
By using Taylor expansion in (12), we have
\begin{equation*}
\tag{13}
\int_{X}\chi(x)(dn+x)^{k}dx-\int_{X}\chi(x)x^{k}dx
=k T_{k-1}(\chi,nd-1),  \quad \hbox{for \ $k, n, d \in \mathbb{N}$}.
\end{equation*}
That is,
\begin{equation*}
B_{k,\chi}(nd)-B_{k,\chi}=kT_{k-1}(\chi,nd-1).
\end{equation*}
Let $w_{1}, w_{2}, d \in \mathbb{N}$. Then we consider the following integral equation
\begin{equation*}
\tag{14}
\begin{split}
&\frac{d\int_{X}\int_{X}\chi(x_{1})\chi(x_{2})e^{(w_{1}x_{1}+
w_{2}x_{2})t}dx_{1}dx_{2}}
{\int_{X}e^{dw_{1} w_{2}xt}dx} \\
&=\frac{t(e^{dw_{1}w_{2}t}-1)}{(e^{w_{1}dt}-1)(e^{w_{2}dt}-1)}
\Big(\sum_{a=0}^{d-1}\chi(a)e^{w_{1}at}\Big)
\Big(\sum_{b=0}^{d-1}\chi(b)e^{w_{2}bt}\Big).
\end{split}
\end{equation*}
From (9) and (12), we note that
\begin{equation*}
\tag{15}
\frac{dw_{1}\int_{X}\chi(x)e^{xt}dx}
{\int_{X}e^{dw_{1}xt}dx}
=\sum_{k=0}^{\infty}\Big(T_{k}(\chi,dw_{1}-1)\Big)\frac{t^{k}}{k!}.
\end{equation*}
Let us consider the $p$-adic functional $T_{\chi}(w_{1},w_{2})$ as follows:
\begin{equation*}
\tag{16}
T_{\chi}(w_{1},w_{2})=
\frac{d\int_{X}\int_{X}\chi(x_{1})\chi(x_{2})e^{(w_{1}x_{1}+w_{2}x_{2}+w_{1}w_{2}x)t}dx_{1}dx_{2}}
{\int_{X}e^{dw_{1}w_{2}x_{3}t}dx_{3}}.
\end{equation*}
Then we see that $T_{\chi}(w_{1},w_{2})$ is symmetric in $w_{1}$ and $w_{2}$, and
\begin{equation*}
\tag{17}
T_{\chi}(w_{1},w_{2})=\frac{t(e^{dw_{1}w_{2}t}-1)e^{w_{1}w_{2}xt}}{(e^{w_{1}dt}-1)(e^{w_{2}dt}-1)}
\Big(\sum_{a=0}^{d-1}\chi(a)e^{w_{1}at}\Big)
\Big(\sum_{b=0}^{d-1}\chi(b)e^{w_{2}bt}\Big).
\end{equation*}
By (16) and (17), we have
\begin{equation*}
\tag{18}
\begin{split}
T_{\chi}(w_{1},w_{2})
&=\Big(\frac{1}{w_{1}}\int_{X}\chi(x_{1})e^{w_{1}(x_{1}+w_{2}x)t} dx_{1}\Big)
\Big(\frac{dw_{1}\int_{X}\chi(x_{2})e^{w_{2}x_{2}t}dx_{2}}
{\int_{X}e^{dw_{1}w_{2}xt}dx}\Big)\\
&=\Big(\frac{1}{w_{1}} \sum_{i=0}^{\infty}B_{i,\chi}(w_{2}x)\frac{w_{1}^{i}t^{i}}{i!} \Big)
\Big( \sum_{k=0}^{\infty}T_{k}(\chi, dw_{1}-1) \frac{w_{2}^{k}t^{k}}{k!} \Big)\\
&=\frac{1}{w_{1}}\Big(\sum_{\ell=0}^{\infty}(\sum_{i=0}^{\ell}
\frac{B_{i,\chi}(w_{2}x)T_{\ell-i}(\chi,dw_{1}-1)w_{1}^{i}w_{2}^{\ell-i}\ell!}{i!(\ell-i)!})
\frac{t^{\ell}}{\ell!}\Big)\\
&=\sum_{\ell=0}^{\infty}\Big(\sum_{i=0}^{\ell}\binom{\ell}{i}B_{i,\chi}(w_{2}x)T_{\ell-i}(\chi,dw_{1}-1)
w_{1}^{i-1}w_{2}^{\ell-i}\Big)\frac{t^{\ell}}{\ell!}.
\end{split}
\end{equation*}
From the symmetric property of $T_{\chi}(w_{1},w_{2})$ in $w_{1}$
and $w_{2}$, we note that
\begin{equation*}
\tag{19}
\begin{split}
T_{\chi}(w_{1},w_{2})
&=\Big(\frac{1}{w_{2}}\int_{X}\chi(x_{2})e^{w_{2}(x_{2}+w_{1}x)t} dx_{2}\Big)
\Big(\frac{dw_{2}\int_{X}\chi(x_{1})e^{w_{1}x_{1}t}dx_{1}}
{\int_{X}e^{dw_{1}w_{2}xt}dx}\Big)\\
&=\Big(\frac{1}{w_{2}} \sum_{i=0}^{\infty}B_{i,\chi}(w_{1}x)\frac{w_{2}^{i}t^{i}}{i!} \Big)
\Big( \sum_{k=0}^{\infty}T_{k}(\chi, dw_{2}-1) \frac{w_{1}^{k}t^{k}}{k!} \Big)\\
&=\frac{1}{w_{2}}\Big(\sum_{\ell=0}^{\infty}(\sum_{i=0}^{\ell}
\frac{B_{i,\chi}(w_{1}x)w_{2}^{i}T_{\ell-i}(\chi,dw_{2}-1)w_{1}^{\ell-i}\ell!}{i!(\ell-i)!})
\frac{t^{\ell}}{\ell!}\Big)\\
&=\sum_{\ell=0}^{\infty}\Big(\sum_{i=0}^{\ell}\binom{\ell}{i}w_{2}^{i-1}w_{1}^{\ell-i}
B_{i,\chi}(w_{1}x)T_{\ell-i}(\chi,dw_{2}-1)
\Big)\frac{t^{\ell}}{\ell!}.
\end{split}
\end{equation*}
By comparing the coefficients on the both sides of (18) and (19), we obtain the following theorem.
{\THM {Theorem {2}}} {\it For $w_{1}, w_{2}, d \in \mathbb{N}$, we have
\begin{equation*}
\begin{split}
&\sum_{i=0}^{\ell}\binom{\ell}{i}B_{i,\chi}(w_{2}x)T_{\ell-i}(\chi,dw_{1}-1)
w_{1}^{i-1}w_{2}^{\ell-i}\\
&=\sum_{i=0}^{\ell}\binom{\ell}{i}
B_{i,\chi}(w_{1}x)T_{\ell-i}(\chi,dw_{2}-1)w_{2}^{i-1}w_{1}^{\ell-i}.
\end{split}
\end{equation*}
\rm}
\par
Let $x=0$ in Theorem 2. Then we have
\begin{equation*}
\begin{split}
&\sum_{i=0}^{\ell}\binom{\ell}{i}B_{i,\chi}T_{\ell-i}(\chi,dw_{1}-1)
w_{1}^{i-1}w_{2}^{\ell-i}\\
&=\sum_{i=0}^{\ell}\binom{\ell}{i}
B_{i,\chi}T_{\ell-i}(\chi,dw_{2}-1)w_{2}^{i-1}w_{1}^{\ell-i}.
\end{split}
\end{equation*}
By (15) and (17), we also see that
\begin{equation*}
\tag{20}
\begin{split}
T_{\chi}(w_{1},w_{2})
&=\Big(\frac{e^{w_{1}w_{2}xt}}{w_{1}}\int_{X}\chi(x_{1})e^{w_{1}x_{1}t} dx_{1}\Big)
\Big(\frac{dw_{1}\int_{X}\chi(x_{2})e^{w_{2}x_{2}t}dx_{2}}
{\int_{X}e^{dw_{1}w_{2}xt}dx} \Big)\\
&=\Big(\frac{e^{w_{1}w_{2}xt}}{w_{1}}\int_{X}\chi(x_{1})e^{w_{1}x_{1}t} dx_{1}\Big)
\Big( \frac{e^{dw_{1}w_{2}t}-1}{e^{w_{2}dt}-1}  \Big)
\Big( \sum_{i=0}^{d-1}\chi(i)e^{w_{2}it}\Big) \\
&=\Big(\frac{e^{w_{1}w_{2}xt}}{w_{1}}\int_{X}\chi(x_{1})e^{w_{1}x_{1}t} dx_{1}\Big)
\Big( \sum_{\ell=0}^{w_{1}-1}\sum_{i=0}^{d-1}e^{w_{2}(i+\ell d)t}\chi(i+\ell d)\Big) \\
&=\Big(\frac{e^{w_{1}w_{2}xt}}{w_{1}}\int_{X}\chi(x_{1})e^{w_{1}x_{1}t} dx_{1}\Big)
\Big(\sum_{i=0}^{dw_{1}-1}e^{w_{2}it}\chi(i)\Big)\\
&=\frac{1}{w_{1}}\sum_{i=0}^{dw_{1}-1}\chi(i)\int_{X}\chi(x_{1})e^{w_{1}(x_{1}+w_{2}x+\frac{w_{2}}{w_{1}}i)t} dx_{1}\\
&=\frac{1}{w_{1}}\sum_{i=0}^{dw_{1}-1}\chi(i)\sum_{k=0}^{\infty}
B_{k,\chi}(w_{2}x+\frac{w_{2}}{w_{1}}i)\frac{w_{1}^{k}t^{k}}{k!}\\
&=\sum_{k=0}^{\infty}\Big(\sum_{i=0}^{dw_{1}-1}\chi(i)B_{k,\chi}(w_{2}x+\frac{w_{2}}{w_{1}}i)w_{1}^{k-1}
\Big)\frac{t^{k}}{k!}.
\end{split}
\end{equation*}
From the symmetric property of $T_{\chi}(w_{1},w_{2})$ in $w_{1}$ and $w_{2}$, we can also derive the following equation.
\begin{equation*}
\tag{21}
\begin{split}
T_{\chi}(w_{1},w_{2})
&=\Big(\frac{e^{w_{1}w_{2}xt}}{w_{2}}\int_{X}\chi(x_{2})e^{w_{2}x_{2}t} dx_{2}\Big)
\Big(\frac{dw_{2}\int_{X}\chi(x_{1})e^{w_{1}x_{1}t}dx_{1}}
{\int_{X}e^{dw_{1}w_{2}xt}dx} \Big)\\
&=\Big(\frac{e^{w_{1}w_{2}xt}}{w_{2}}\int_{X}\chi(x_{2})e^{w_{2}x_{2}t} dx_{2}\Big)
\Big( \frac{e^{dw_{1}w_{2}t}-1}{e^{w_{1}dt}-1}  \Big)
\Big( \sum_{i=0}^{d-1}\chi(i)e^{w_{1}it}\Big) \\
&=\Big(\frac{e^{w_{1}w_{2}xt}}{w_{2}}\int_{X}\chi(x_{2})e^{w_{2}x_{2}t} dx_{2}\Big)
\Big( \sum_{\ell=0}^{w_{2}-1}e^{w_{1}d\ell t}\Big)
\Big( \sum_{i=0}^{d-1}\chi(i)e^{w_{1}it}\Big) \\
&=\frac{1}{w_{2}}\sum_{i=0}^{dw_{2}-1}\chi(i)\int_{X}\chi(x_{2})e^{w_{2}(x_{2}+w_{1}x+\frac{w_{1}}{w_{2}}i)t} dx_{2}\\
&=\frac{1}{w_{2}}\sum_{i=0}^{dw_{2}-1}\chi(i)\sum_{k=0}^{\infty}B_{k,\chi}(w_{1}x+\frac{w_{1}}{w_{2}}i)
\frac{w_{2}^{k}t^{k}}{k!}\\
&=\sum_{k=0}^{\infty}\Big\{ \sum_{i=0}^{dw_{2}-1}\chi(i)B_{k,\chi}(w_{1}x+\frac{w_{1}}{w_{2}}i)w_{2}^{k-1}\Big\}
\frac{t^{k}}{k!}.
\end{split}
\end{equation*}
By comparing the coefficients on the both sides of (20) and (21), we obtain the following theorem.
{\THM {Theorem {3}}} {\it For $w_{1}, w_{2}, d \in \mathbb{N}$, we have
\begin{equation*}
\sum_{i=0}^{dw_{1}-1}\chi(i)B_{k,\chi}(w_{2}x+\frac{w_{2}}{w_{1}}i)w_{1}^{k-1}
=\sum_{i=0}^{dw_{2}-1}\chi(i)B_{k,\chi}(w_{1}x+\frac{w_{1}}{w_{2}}i)w_{2}^{k-1}.
\end{equation*}
\rm}
\medskip
{\REM {Remark}} {\it Let $x=0$ in ${\rm Theorem\ 3}$. Then we see that
\begin{equation*}
\sum_{i=0}^{dw_{1}-1}\chi(i)B_{k,\chi}(\frac{w_{2}}{w_{1}}i)w_{1}^{k-1}
=\sum_{i=0}^{dw_{2}-1}\chi(i)B_{k,\chi}(\frac{w_{1}}{w_{2}}i)w_{2}^{k-1}.
\end{equation*}
\par
If we take $w_{2}=1$, then we have
\begin{equation*}
\sum_{i=0}^{dw_{1}-1}\chi(i)B_{k,\chi}(\frac{i}{w_{1}})w_{1}^{k-1}
=\sum_{i=0}^{d-1}\chi(i)B_{k,\chi}({w_{1}}i).
\end{equation*}
\rm}
\bigskip\bigskip
\begin{center}\begin{large}
{\sc References}
\end{large}\end{center}
\par
\begin{enumerate}
\item[{[1]}] L. C. Carlitz,{\it $q$-Bernoulli numbers and polynomials}, Duke Math. J.
{\bf 15} (1948), 987-1000.

\item[{[2]}] M. Cenkci, Y. Sisek, V. Kurt, {\it Further remarks on multiple $p$-adic $q$-$L$-function of two variables}, Adv. Stud. Contemp. Math. {\bf 14} (2007), 49-68.

\item[{[3]}] M. Cenkci, M. Can, V. Kurt, {\it Multiple two-variable $q$-$L$-function and its behavior at $s=0$}, Russ. J. Math. Phys. {\bf 15(4)} (2008), 447-459.

\item[{[4]}] T. Ernst, {\it Example of a $q$-umbral calculus}, Adv. Stud. Contemp. Math. {\bf 16(1)} (2008), 1-22.

\item[{[5]}] A. S. Hegazi, M. Mansour, {\it A note on $q$-Bernoulli numbers and polynomials}, J. Nonlinear Math. Phys. {\bf 13} (2006), 9-18.

\item[{[6]}] T. Kim, {\it $q$-Volkenborn Integration}, Russ. J. Math. Phys. {\bf 9} (2002), 288-299.

\item[{[7]}] T. Kim, {\it Non-archimedean $q$-integrals associated with multiple Changhee $q$-Bernoulli polynomials }, Russ. J. Math. Phys. {\bf 10} (2003), 91-98.

\item[{[8]}] T. Kim, {\it Power series and asymptotic series associated with the $q$-analog of the two-variable
$p$-adic $L$-function}, Russ. J. Math. Phys. {\bf 12(2)} (2005), 186-196.

\item[{[9]}] T. Kim, {\it Multiple $p$-adic $L$-function} Russ. J. Math. Phys. {\bf 13(2)} (2006), 151-157.

\item[{[10]}] T. Kim, {\it $q$-Euler numbers and polynomials associated with $p$-adic $q$-integral}, J. Nonlinear Math. Phys.{\bf 14 (1)} (2007), 15-27.

\item[{[11]}] T. Kim, {\it A note on $p$-adic $q$-integral on $\mathbb{Z}_{p}$}, Adv. Stud. Contemp. Math. {\bf 15} (2007), 133-138.

\item[{[12]}] T. Kim, {\it $q$-Bernoulli numbers and polynomials associated with Gaussian binomial coefficients}, Russ. J. Math. Phys. {\bf 15 (1)} (2008), 51-57.

\item[{[13]}] T. Kim, {\it On the symmetry of the $q$-Bernoulli polynomials}, Abstr. Appl. Anal. {\bf 2008} (2008), Article ID914367, 7 pages.

\item[{[14]}] T. Kim, {\it Symmetries $p$-adic invariant integral on $\mathbb{Z}_{p}$ for Bernoulli and Euler polynomials}, J. Difference Equ. Appl. {\bf 14 (12)} (2008), 1267-1277.

\item[{[15]}] T. Kim, {\it Note on $q$-Genocchi numbers and polynomials}, Adv. Stud. Contemp. Math. {\bf 17 (1)} (2008), 9-15.

\item[{[16]}] T. Kim, {\it Symmetry of power sum polynomials and multivariate fermionic $p$-adic invariant integral on $\mathbb{Z}_{p}$}, Russ. J. Math. Phys. {\bf 16 (1)} (2009), 51-54.

\item[{[17]}] Y. -H. Kim, W. Kim, L.-C. Jang, {\it On the $q$-extension of Apostol-Euler numbers and polynomials}, Abstr. Appl. Anal. {\bf 2008} (2008), Article ID296159, 10 pages.

 \item[{[18]}] B. A. Kupershmidt, {\it Reflection symmetries of $q$-Bernoulli polynomials}, J. Nonlinear Math. Phys. {\bf 12} (2005), 412-422.

\item[{[19]}] H. Ozden, Y. Simsek, S. -H. Rim, I. N. Cangul {\it A note on $p$-adic $q$-Euler measure}, Adv. Stud. Contemp. Math. {\bf 14} (2007), 233-239.

\item[{[20]}] K. H. Park, Y.-H. Kim {\it On some arithmetical properties of the Genocchi numbers and polynomials}, Advances in Difference Equations. http://www.hindawi.com/journals/
\newline
ade/aip.195049.html.

\item[{[21]}] M. Schork, {\it A representation of the $q$-fermionic commutation relations and the limit $q=1$},
Russ. J. Math. Phys. {\bf 12 (3)} (2005), 394-399.

\item[{[22]}] Y. Simsek, {\it Theorems on twisted $L$-function and twisted Bernoulli numbers },
Adv. Stud. Contemp. Math. {\bf 11} (2005), 205-218.

\item[{[23]}] Y. Simsek, {\it On $p$-adic twisted $q$-$L$-functions related to generalized twisted Bernoulli numbers}, Russ. J. Math. Phys. {\bf 13(3)} (2006), 340-348.

\item[{[24]}] Y. Simsek, {\it Complete sums of $(h,q)$-extension of the Euler polynomials and numbers}, arXiv:0707.2849vl[math.NT].

\item[{[25]}] Y.-H. Kim, K.-W. Hwang. {\it Symmetry of power sum and twisted  Bernoulli polynomials},
Adv. Stud. Contemp. Math. {\bf 18 (2)} (2009), 105-113.

\end{enumerate}

\newpage
\par\noindent
\mpn { \bpn {\small Taekyun {\sc Kim} \mpn
Division of General Education-Mathematics,
\pn Kwangwoon University, Seoul 139-701, S. Korea
\pn {\it
E-mail:}\ {\sf tkkim@kw.ac.kr} }

\end{document}